\documentstyle[amssymb,12pt]{amsart}

\newtheorem{cor}{Corollary}
\newtheorem{lem}{Lemma}
\newtheorem{prop}{Proposition}

\theoremstyle{definition}

\theoremstyle{definition}
\newtheorem{thm}{Theorem}
\newtheorem{conj}{Conjecture}

\theoremstyle{remark}
\newtheorem{rem}{Remark}

\numberwithin{equation}{section}

\begin{document}

\newcommand{\thmref}[1]{Theorem~\ref{#1}}
\newcommand{\secref}[1]{Sect.~\ref{#1}}
\newcommand{\lemref}[1]{Lemma~\ref{#1}}
\newcommand{\propref}[1]{Proposition~\ref{#1}}
\newcommand{\corref}[1]{Corollary~\ref{#1}}
\newcommand{\remref}[1]{Remark~\ref{#1}}
\newcommand{\nc}{\newcommand}
\nc{\on}{\operatorname}
\nc{\ch}{\mbox{ch}}
\nc{\Z}{{\Bbb Z}}
\nc{\C}{{\Bbb C}}
\nc{\pone}{{\Bbb C}{\Bbb P}^1}
\nc{\pa}{\partial}
\nc{\F}{{\cal F}}
\nc{\arr}{\rightarrow}
\nc{\larr}{\longrightarrow}
\nc{\al}{\alpha}
\nc{\ri}{\rangle}
\nc{\lef}{\langle}
\nc{\W}{{\cal W}}
\nc{\la}{\lambda}
\nc{\ep}{\epsilon}
\nc{\su}{\widehat{\goth{sl}}_2}
\nc{\sw}{\goth{sl}}
\nc{\g}{\goth{g}}
\nc{\h}{\goth{h}}
\nc{\n}{\goth{n}}
\nc{\N}{\widehat{\n}}
\nc{\ab}{\goth{a}}
\nc{\G}{\widehat{\g}}
\nc{\De}{\Delta_+}
\nc{\gt}{\widetilde{\g}}
\nc{\Ga}{\Gamma}
\nc{\one}{{\bold 1}}
\nc{\hh}{\widehat{\h}}
\nc{\z}{{\goth Z}}
\nc{\zz}{{\cal Z}}
\nc{\Hh}{{\cal H}}
\nc{\qp}{q^{\frac{k}{2}}}
\nc{\qm}{q^{-\frac{k}{2}}}
\nc{\La}{\Lambda}
\nc{\wt}{\widetilde}
\nc{\qn}{\frac{[m]_q^2}{[2m]_q}}
\nc{\cri}{_{\on{cr}}}
\nc{\k}{h^\vee}
\nc{\sun}{\widehat{\sw}_N}
\nc{\HH}{{\cal H}_h(\sw_N)}
\nc{\ca}{\wt{{\cal A}}_{h,k}(\sw_2)}
\nc{\si}{\sigma}
\nc{\gl}{\widehat{\goth{g}\goth{l}}_2}
\nc{\el}{\ell}
\nc{\s}{s}
\nc{\bi}{\bibitem}

\title[Quantum algebras and ${\cal W}$--algebras]{Quantum affine algebras
and deformations of the Virasoro and ${\cal W}$--algebras}

\author{Edward Frenkel}\thanks{Partially supported by NSF grants
DMS-9205303 and DMS-9296120}
\address{Department of Mathematics, Harvard University, Cambridge, MA
02138, USA}

\author{Nikolai Reshetikhin}
\address{Department of Mathematics, University of California, Berkeley, CA
94720, USA}

\date{May 1995}

\maketitle

\section{Introduction.}
\subsection{} In this paper we generalize some results concerning affine
Kac-Moody algebras at the critical level to the corresponding quantized
universal enveloping algebras. Here is the short description of these
results for the affine algebras.

\renewcommand{\labelenumi}{\normalshape(\roman{enumi})}
\begin {enumerate}
\item
Let $\wt{U}(\G)\cri$ be a completion of the universal enveloping algebra
of an affine algebra $\G$ at the critical level $-h^\vee$ (the precise
definition is given in \S~2). This algebra possesses a large center
$Z(\G)$, which has a natural Poisson structure. B.~Feigin and E.~Frenkel
have shown that $Z(\G)$ is isomorphic to the classical $\W$--algebra
$\W(\g^L)$ associated to the simple Lie algebra $\g^L$, which is Langlands
dual to $\g$ \cite{FF:ctr}.

\item
The $\W$--algebra $\W(\g^L)$ consists of functionals on a certain Poisson
manifold ${\cal C}(\g^L)$ obtained by the Drinfeld-Sokolov hamiltonian
reduction \cite{DS} from a hyperplane in the dual space to the affine
algebra $\widehat{\g^L}$. Elements of ${\cal C}(\g^L)$, called
$\g^L$--opers in \cite{BD}, can be considered as connections on a certain
$G^L$--bundle over the circle with some extra structure. To a $\g^L$--oper
one can attach a $\G$--module, on which the center acts according to the
corresponding character. These $\G$--modules can be considered as analogues
of admissible representations of a simple group over a local
non-archimedian field. They can be used in carrying out the geometric
Langlands correspondence proposed by A.~Beilinson and V.~Drinfeld
\cite{BD}.

\item
The Wakimoto realization of $\G$ \cite{Wak,FF:ff,FF:ctr,F:talk} provides a
map from $\wt{U}(\G)\cri$ to the tensor product of a certain Heisenberg
algebra and a commutative algebra ${\cal H}(\g)$. The restriction of this
map to $Z(\G)$ gives us a homomorphism $Z(\G) \arr {\cal H}(\g)$, which is
an analogue of the Harish-Chandra homomorphism. The corresponding map
$\W(\g^L) \arr {\cal H}(\g)$ is nothing but the Miura transformation, which
has been defined for an arbitrary $\g$ by V.~Drinfeld and V.~Sokolov
\cite{DS}.

\item
The algebra ${\cal H}(\g)$ consists of functionals on a hyperplane ${\cal
F}$ in the dual space to the homogeneous Heisenberg subalgebra $\hh$ of
$\G$. The algebra ${\cal H}(\g)$ is the classical limit of a completion
of $U(\hh)$, and hence it is a Heisenberg-Poisson algebra. The Miura
transformation $\W(\g^L) \arr {\cal H}(\g)$ is a homomorphism of Poisson
algebras.

\end{enumerate}

For example, the center of $\wt{U}(\su)\cri$ is generated by the Sugawara
operators, and ${\cal C}(\sw_2)$ is isomorphic to a hyperplane ${\cal L}$
in the dual space to the Virasoro algebra. Thus, the center of
$\wt{U}(\su)\cri$ is isomorphic to the Poisson algebra $\W(\sw_2)$ of
functionals on ${\cal L}$. The Poisson structure on ${\cal L}$ is often
called the second Gelfand-Dickey structure. We call $\W(\sw_2)$ the
classical Virasoro algebra.

Elements of ${\cal L}$ can be considered as projective connections on the
circle, i.e. differential operators of the form $\pa_t^2 - q(t)$; these are
the $\sw_2$--opers. On the other hand, elements of ${\cal F}$ can be
considered as connections on a rank one bundle over the circle,
i.e. differential operators of the form $\pa_t - \frac{1}{2}\chi(t)$. The
Miura transformation sends a connection $\pa_t - \frac{1}{2}\chi(t)$ to the
projective connection
\begin{equation}    \label{miura}
\pa_t^2 - q(t) = (\pa_t - \frac{1}{2}\chi(t))(\pa_t + \frac{1}{2}\chi(t)).
\end{equation}
This gives us a homomorphism of Poisson algebras $\W(\sw_2) \arr {\cal
H}(\sw_2)$.

In \cite{FFR,Fr} these results were used to give a new interpretation of
the Bethe ansatz in the Gaudin models of statistical mechanics. This
allowed to gain new insights on completeness of Bethe ansatz, and to relate
Bethe ansatz to the geometric Langlands correspondence.

\subsection{} There are many indications that these results can be
generalized to the center of a completion $\wt{U}_h(\G)\cri$ of the
quantum affine algebra $U_q(\G)$ at the critical level. An
explicit construction of central elements of a quantum affine algebra at
the critical level has been given by N.~Reshetikhin and
M.~Semenov-Tian-Shansky \cite{RS}. Later, J.~Ding and P.~Etingof \cite{DE}
showed that those elements generate all singular vectors of imaginary
weight in Verma modules over $U_h(\G)\cri$. This makes us to believe that
the center of $\wt{U}_h(\G)\cri$ is generated, in an appropriate sense, by
the elements constructed in \cite{RS}.

The center $Z_h(\G)$ of $\wt{U}_h(\G)\cri$ possesses a natural Poisson
structure, which is a $q$--deforma\-tion of the Poisson structure on
$Z(\G)$. A natural question is to describe $Z_h(\G)$ and its spectrum.

In this paper we do this explicitly for $\wt{U}_h(\su)\cri$ by using its
Wakimoto realization. Our results for $\wt{U}_h(\su)\cri$ can be summarized
as follows.

\renewcommand{\labelenumi}{\normalshape(\roman{enumi})}
\begin {enumerate}
\item
The center $Z_h(\su)$ of $\wt{U}_h(\su)\cri$ contains the Fourier coefficients
of a power series $\el(z)$ given in \cite{RS} in terms of the
Reshetikhin-Semenov-Tian-Shansky (RS) realization of $U_h(\su)\cri$. We
rewrite $\el(z)$ in terms of the Drinfeld realization \cite{Dr:new}, using
the explicit isomorphism between the two realizations established by
J.~Ding and I.~Frenkel \cite{DF}. This gives us a $q$--analogue of the
Sugawara formula.

\item
Wakimoto realizations of $U_q(\su)$ have been given in \cite{qWak} in terms
of the Drinfeld realization; we use the presentation due to H.~Awata,
S.~Odake, and J.~Shiraishi \cite{AOS}. It gives us a homomorphism from
$\wt{U}_h(\G)\cri$ to the tensor product of a certain Heisenberg algebra
and a Heisenberg-Poisson algebra ${\cal H}_h(\sw_2)$. Its restriction to
$Z_h(\su)$ is a homomorphism of Poisson algebras $Z_h(\su) \arr {\cal
H}_h(\sw_2)$, which is a $q$--deformation of the Miura transformation.

\item
We find the image of $\el(z)$ in ${\cal H}_h(\sw_2)$:
\begin{equation}    \label{map}
\el(z) \longrightarrow \s(z) = \Lambda(zq) + \Lambda(zq^{-1})^{-1},
\end{equation}
where $\La(z)$ is a generating function of elements of ${\cal
H}_h(\sw_2)$. Using the Poisson structure on ${\cal H}_h(\sw_2)$ we
compute the Poisson structure on $Z_h(\su)$:
\begin{multline}    \label{mult}
\{ \el(z),\el(w) \} = \\ 2h (q-q^{-1}) \sum_{m\in\Z} \left( \frac{w}{z}
\right)^m \qn \; \el(z) \el(w) + 2h \sum_{m\in\Z} \left( \frac{w}{zq^2}
\right)^m - 2h \sum_{m\in\Z} \left( \frac{wq^2}{z} \right)^m.
\end{multline}
This gives us a $q$--deformation of the classical Virasoro algebra.
\end{enumerate}

Formula \eqref{map} shows that if we attach to $\La(z)$ a first order
$q$--difference operator ${\cal D}_q - \La(zq)$, where $[{\cal D}_q \cdot
f](z) = f(zq^2)$, then to $\s(z)$ we can attach in a natural way a second
order $q$--difference operator of the form ${\cal D}_q + {\cal D}_q^{-1} -
\s(z)$. Indeed, let $Q(z)$ be a solution of the difference equation $Q(zq)
= \La(z) Q(zq^{-1})$. Then from formula \eqref{map} we obtain $$({\cal D}_q
+ {\cal D}_q^{-1} - \s(z)) Q(z) = 0.$$ The latter equation written as
\begin{equation}    \label{baxter}
\s(z) = \frac{Q(zq^2)}{Q(z)} + \frac{Q(zq^{-2})}{Q(z)}
\end{equation}
was used by R.~Baxter \cite{Baxter} in his study of the eight vertex
model. Similar formulas appeared in \cite{Lieb} as the result of
computation of the spectrum of the transfer-matrix of the six vertex model,
an integrable model of statistical mechanics associated to $U_q(\su)$. In
this context, the function $Q(z)$ is a product of a ``vacuum value'' and a
polynomial, whose zeroes are solutions of Bethe ansatz equations.

\subsection{} Thus, we have interpreted formulas \eqref{map} and
\eqref{baxter} as a hamiltonian map, which can be considered as a
$q$--analogue of the Miura transformation. In fact, the Miura
transformation plays the same role as Baxter's formula \eqref{baxter} in
the Gaudin models, cf. \cite{S,FFR,Fr}.

The Miura transformation \eqref{miura} is the classical limit of the free
field realization of the Virasoro algebra. Free field realizations play an
important role in conformal field theory, cf. \cite{F:talk}. It is quite
remarkable that a $q$--analogue of free field realization appears in the
context of Bethe ansatz in statistical mechanics.

Analogues of formula \eqref{map} for transfer-matrices of integrable models
associated to other quantum affine algebras are known,
cf. e.g. \cite{KuR,BR} for $U_q(\sun)$, \cite{R1,R2,KS} for other
$U_q(\G)$. Motivated by our computation in the case of $U_q(\su)$ we expect
that the formulas for the $q$--deformation of the Miura transformation of
the center of $\wt{U}_h(\G)\cri$ coincide with the formulas for the
transfer-matrices corresponding to $U_q(\G)$.

In particular, for $\g=\sw_N$ we obtain the following picture.

\renewcommand{\labelenumi}{\normalshape(\roman{enumi})}
\begin {enumerate}
\item
In \cite{RS} the generating functions of central elements
$\el_1(z),\ldots,\el_{N-1}(z)$ of $\wt{U}_h(\sun)\cri$ corresponding to the
fundamental representations have been constructed. The Fourier coefficients
of $\el_i(z)$'s generate a central subalgebra $Z_h(\sun)$ of
$\wt{U}_h(\sun)\cri$, which is closed with respect to the Poisson structure.

\item
The Wakimoto realization of $U_h(\sun)\cri$ \cite{AOS} gives rise to a
homomorphism of Poisson algebras $Z_h(\sun) \arr {\cal H}_h(\sw_N)$, where
${\cal H}_h(\sw_N)$ is a Heisenberg-Poisson algebra. This is a
$q$--deformation of the Miura transformation. We find a formula for the
image $\s_i(z)$ of each generating function $\el_i(z)$ in ${\cal
H}_h(\sw_N)$. These formulas match formulas for the spectra of the
corresponding transfer-matrices in integrable models associated to
$U_q(\sun)$ \cite{KuR,BR}.

\item
We explicitly compute the Poisson brackets between $\s_i(z)$'s
in ${\cal H}_h(\sw_N)$ generalizing formula \eqref{mult}. Thus, we obtain
an interesting Poisson subalgebra $\W_h(\sw_N)$ of the Heisenberg-Poisson
algebra ${\cal H}_h(\sw_N)$, which is a $q$--deformation of the classical
$\W$--algebra $\W(\sw_N)$.

\item
Recall that elements of the spectrum of $\W(\sw_N)$ can be considered as
$N$th order differential operators. We show that elements of the spectrum
of $\W_h(\sw_N)$ can be considered as $N$th order $q$--difference operators
of the form $${\cal D}_q^N - \s_{N-1}(z) {\cal D}_q^{N-1} + \s_{N-2}(z) {\cal
D}_q^{N-2} - \ldots - (-1)^N \s_1 {\cal D}_q + (-1)^N.$$
\end{enumerate}

We generalize (i) and (ii) to all quantum affine algebras of classical
types. The computation of Poisson brackets is straightforward, and will be
given in our next paper \cite{next} along with results regarding quantum
affine algebras of exceptional types.

Using our results in the same way as in \cite{FFR} we can give a new
interpretation of the Bethe ansatz in integrable models associated to
quantum affine algebras. This and other applications will be discussed in
\cite{next}.

The paper is organized as follows. In \S~2 we recall results concerning the
center of $\wt{U}(\su)\cri$ and Miura transformation. In \S\S~3-5 we
consider the Drinfeld and the RS realizations of $U_h(\su)_k$ and the
isomorphism between them. In \S~6 we rewrite the RS formula for the
generating function $\el(z)$ of central elements in $\wt{U}_h(\su)\cri$ in
terms of the Drinfeld realization. In \S\S~7-9 we recall the Wakimoto
realization of $U_h(\su)_k$, and use it to find an explicit formula for the
image of $\el(z)$ in ${\cal H}_h(\sw_2)$ and to compute the Poisson bracket
on $Z_h(\su)$. In \S\S~10 and 11 we generalize these results to
$U_h(\sun)\cri$ and other quantum affine algebras of classical types.

\section{The center of $\wt{U}(\su)$ at the critical level.}
\subsection{The structure of the center.}
Let $\g$ be a finite-dimensional simple Lie algebra. The affine algebra
$\G$ is the extension of $\g \otimes \C[t,t^{-1}]$ by a one-dimensional
center $\C K$. For $A \in \g, n\in\Z$, denote $A[n] = A \otimes t^n$ and
put $$A(z) = \sum_{n\in\Z} A[n] z^{-n-1}.$$

Introduce a completion $\wt{U}(\G)$ of $U(\G)$, the universal enveloping
algebra of $\G$: $$\wt{U}(\G) = \underset{\longleftarrow}{\lim} \;
U(\G)/U(\G) (\g \otimes t^n\C[t]), \quad \quad n>0.$$ This is an
associative algebra. It acts on $\G$--modules $M$ which satisfy the
following property: for any $x \in M$ there exists $N\in\Z_+$ such that
$A[n] \cdot x = 0$ for any $A \in \g$ if $n>N$. For $k\in\C$ put
$\wt{U}(\G)_k=\wt{U}(\G)/(K-k)$, and let $\wt{U}(\G)\cri$ be
$\wt{U}(\G)_{-\k}$.

The algebra $\wt{U}(\G)_k$ contains the local completion
$U(\G)_{k,\on{loc}}$ introduced in \cite{FF:ctr}. The center of
$U(\G)_{k,\on{loc}}$ has been described in \cite{FF:ctr}. It consists of
the constants when $k\neq -\k$, but becomes ``large'' when $k=-\k$. Let us
recall its description in the case $\G=\su$; in this case $\k=2$.

Let $\{ e, h, f \}$ be the standard basis of $\sw_2$. Introduce the
generating function of the Sugawara operators $S_n$ by formula
$$S(z) = \sum_{n \in \Z} S_n z^{-n-2} = \frac{1}{4} :h(z)^2: + \frac{1}{2}
:e(z) f(z): + \frac{1}{2} :f(z) e(z):.$$ It is well-known that $$[S_n,A[m]]
= -(k+2) m A[n+m]$$ for any $A\in\g$, and
\begin{equation}    \label{snsm}
[S_n,S_m] = (k+2) \left[ (n-m) S_{n+m} + \frac{k}{4} (n^3-n) \delta_{n,-m}
\right].
\end{equation}
Therefore, if $k
\neq -2$, the operators $L_n = S_n/(k+2)$ generate the Virasoro algebra. If
$k=-2$, the operators $S_n, n \in \Z$, are central elements of
$U(\su)_{-\k,\on{loc}}$.

There is a natural Poisson structure on the center $Z(\G)$ of
$\wt{U}(\G)\cri$: for any $A, B \in Z(\G)$, let $A', B'$ be their liftings
to $\wt{U}(\su)$. Then we have $[A',B'] = (K+\k) C' +
(K+\k)^2(\ldots)$. Let $C$ be the projection of $C' \in \wt{U}(\su)$ to
$\wt{U}(\G)\cri$. Then the formula $\{ A,B \} = C$ defines a Poisson
bracket on $Z(\G)$, which does not depend on the liftings.

We obtain from formula \eqref{snsm}:
\begin{equation}     \label{sugpb}
\{ S_n,S_m \} = (n-m) S_{n+m} - \frac{1}{2} (n^3-n) \delta_{n,-m}.
\end{equation}

Consider the hyperplane ${\cal L}$ in the dual space to the Virasoro
algebra, which consists of those linear functionals on the Virasoro algebra
which take value $-6$ on the central element (this corresponds to the
factor $-1/2=-6/12$ in the second term of formula \eqref{sugpb}). This
hyperplane is isomorphic to the space of projective connections on the
formal punctured disc $\pa_z^2 - q(z)$, where $q(z) \in \C((z))$, in the
sense that the natural action of vector fields on it coincides with the
coadjoint action of the Virasoro algebra on ${\cal L}$. Let $\W(\sw_2)$ be
the Poisson algebra of local functionals on ${\cal L}$. It is the classical
limit of the local completion of the universal enveloping algebra of the
Virasoro algebra. Therefore we call $\W(\sw_2)$ the classical Virasoro
algebra.

In the case $\G=\su$ the result of \cite{FF:ctr} is that the center $Z(\su)$ of
$U(\su)_{-2,\on{loc}}$ is isomorphic to $\W(\sw_2)$. This isomorphism sends
$S_n$ to the local functional $\pa_z^2 - q(z) \arr \int q(z) z^{n+1}
dz$. According to formula \eqref{sugpb}, this is an isomorphism of Poisson
algebras.

\subsection{Wakimoto modules and Miura transformation.}
Consider the Heisenberg algebra $\Gamma$ with generators $a_n, a^*_n,
n\in\Z$, and relations $$[a_n,a_m] = [a^*_n,a^*_m] = 0, \quad \quad
[a_n,a^*_m] = \delta_{n,-m},$$ and the Heisenberg algebra ${\cal
H}'(\sw_2)$ with generators $\chi_n, n\in\Z$, and ${\bold 1}$ and relations
$$[\chi_n,\chi_m] = 2 n \delta_{n,-m} \one, \quad \quad [\chi_n,\one] =
0.$$ Introduce the generating functions $$a(z) = \sum_{n\in\Z} a_n
z^{-n-1}, \quad \quad a^*(z) = \sum_{n\in\Z} a^*_n z^{-n} \quad \quad
\chi(z) = \sum_{n\in\Z} \chi_n z^{-n-1}.$$ We define an embedding $\phi$ of
$\su$ into a completion of $\Gamma \otimes {\cal H}'(\sw_2)$ by the
formulas:
\begin{gather*}
\phi[e(z)] = a(z),\\ \phi[h(z)] = -2 :a(z) a^*(z): + \chi(z),\\ \phi[f(z)]
= - :a(z) a^*(z) a^*(z): + ({\bold 1} - 2) \pa_z a^*(z) + \chi(z) a^*(z),\\
\phi(K) = {\bold 1} - 2.
\end{gather*}

The algebra $\Gamma$ has a standard Fock representation $M$ generated by a
vector $v$, such that $$a_n v = 0, \quad n\geq 0; \quad \quad a^*_n v = 0,
\quad n>0.$$

The algebra ${\cal H}'(\sw_2)$ has a family of Fock representations
$\pi_{\mu,\kappa}, \mu \in \C, \kappa \in \C, \kappa \neq 0$, generated by
a vector $v_{\mu,\kappa}$, such that
$$\chi_n v_{\mu,\kappa} = \mu \delta_{n,0} v_{\mu,\kappa}, \quad n\geq 0;
\quad \quad {\bold 1} v_{\mu,\kappa} = \kappa v_{\mu,\kappa}.$$

It also has an infinite-dimensional family of one-dimensional
representations $\C_{x(z)}$, where $x(z) = \sum_{n\in\Z} x_n z^{-n-1} \in
\C((z))$, on which $\chi_n$ acts by multiplication by $x_n$, and ${\bold
1}$ acts by $0$. Using the homomorphism $\phi$ we obtain representations of
$\su$ of non-critical level in $M \otimes \pi_{\mu,\kappa}$, and
representations of critical level in $M\otimes \C_{x(z)}$. These are called
Wakimoto modules.

The algebra ${\cal H}'(\sw_2)$ can be considered as a deformation of the
commutative algebra $\Hh(\sw_2) = {\cal H}'(\sw_2)/({\bold 1})$. It induces a
Poisson structure on $\Hh(\sw_2)$, which is called the Heisenberg-Poisson
structure. It is determined by the following Poisson brackets of the
generators:
\begin{equation}    \label{pb}
\{ \chi_n,\chi_m \} = 2n \delta_{n,-m}.
\end{equation}

The homomorphism $\phi$ defines a homomorphism $\phi_{-2}$ from $\su$ to a
completion of $\Gamma \otimes \Hh(\sw_2)$, which maps $K$ to $-2$.  One can
check that under $\phi_{-2}$, the Sugawara series $S(z)$ is mapped to
\begin{equation}    \label{mm}
\frac{1}{4} \chi(z)^2 - \frac{1}{2} \pa_z \chi(z).
\end{equation}
Therefore $\phi_{-2}$ defines an embedding of $Z(\su)$ into a completion of
${\cal H}(\sw_2)$, which we call the Miura transformation.

The Poisson structure between central elements $A, B \in Z(\su)$ has been
defined via the commutator of their liftings to $\wt{U}(\su)$. The result
does not depend on the choice of a lifting. Moreover, we will obtain the
same result if we take the commutator between the liftings of the images
$A', B'$ of $A, B$ in the completion of $\Gamma \otimes {\cal
H}(\sw_2)$. Since the image of $Z(\su)$ actually lies in the completion of
$\Hh(\sw_2)$, we can take liftings lying in the completion of ${\cal
H}'(\sw_2)$. But then the image of the Poisson bracket between $A, B \in
Z(\su)$ in $\Hh(\sw_2)$ will coincide with the Poisson bracket between $A',
B' \in \Hh(\sw_2)$. Hence the Miura transformation $Z(\su) \arr \Hh(\sw_2)$
is a homomorphism of Poisson algebras.

Therefore we can compute the Poisson structure on $Z(\su)$ using formulas
\eqref{mm} and \eqref{pb}. This gives us the Poisson bracket between the
Sugawara operators, which coincides with formula \eqref{sugpb}.

\section{The Drinfeld realization.}
Let $U_q(\su)$ be the associative algebra over the ring $\C[q,q^{-1}]$,
with generators $e_i, f_i$, and $K_i, K_i^{-1}, i=0,1$, which satisfy the
following relations \cite{D,J}:

$$K_i K_j=K_j K_i,$$ $$K_i e_j=q^{A_{ij}} e_j K_i, \quad \quad
K_if_j=q^{-A_{ij}}f_jK_i,$$
$$[e_i,f_j] = {\delta}_{ij}(q-q^{-1})(K_i-K_i^{-1}),$$
$$e_i^3e_j-(q+1+q^{-1})(e_i^2e_je_i-e_ie_je_i^2)-e_je_i^3=0,$$
$$f_i^3f_j-(q+1+q^{-1})(f_i^2f_jf_i-f_if_jf_i^2)-f_jf_i^3=0,$$ where
$A_{ij}, i,j=0,1$, are the entries of the Cartan matrix of $\su$:
$A_{00}=A_{11}=-A_{01}=-A_{10}=2$.

For any $h, k\in\C$ consider the quotient $U_h(\su)_k$ of $U_q(\su)$ by the
relations $q=e^h, K_0 K_1 = e^{hk}$. This is the Drinfeld-Jimbo realization
of $U_h(\su)_k$.

There exist two other realizations of $U_h(\su)_k$: Drinfeld's realization
\cite{Dr:new}, and Resheti\-khin-Semenov-Tian-Shansky (RS) realization. The
equivalence between the Drinfeld-Jimbo and the Drinfeld realizations has
been established by V.~Drinfeld \cite{Dr:new}, cf. also
\cite{KhT,B,LSS}. The equivalence between the Drinfeld and the RS
realizations has been established by J.~Ding and I.~Frenkel \cite{DF}. The
equivalence between the Drinfeld-Jimbo and the RS realization follows from
these two equivalences, but it can also be established directly along the
lines of \cite{DF}.

First we consider the Drinfeld realization. It is important for us because
the Wakimoto realization is defined in terms of this realization.

Introduce formal power series in $x$
\begin{equation}    \label{f}
f(x) = \frac{(x;q^4)(xq^4;q^4)}{(xq^2;q^4)^2},
\end{equation}
where
$$(a;b) = \prod_{n=0}^\infty (1-ab^n).$$

\begin{rem}    \label{power}
Each coefficient of $f(x)$ is itself a series in $q$, which converges for
$|q|<1$ and can be analytically continued to the whole complex plane except
for some roots of unity. Thus, we can extend $f(x)$ as formal power series
in $x$ to all $q$ except for the roots of unity. In what follows we will
exclude roots of unity from consideration.  \qed\end{rem}

Let $h \in \C \backslash \{ 2\pi i {\Bbb Q} \}, k \in \C$; put $q=e^h$. We
define an associative algebra $U_h(\gl)_k$ over $\C$ with generators $E[n],
F[n], n\in\Z$, and $k^\pm_i[n], i=1,2; n\in \mp\Z_+$. Introduce the
generating functions $$E(z) = \sum_{n\in\Z} E[n] z^{-n}, \quad F(z) =
\sum_{n\in\Z} F[n] z^{-n}, \quad k^\pm_i(z) = \sum_{n=0}^\infty k^\pm_i[\mp
n] z^{\pm n}.$$ The defining relations in $U_h(\gl)_k$ are
$$k^+_i[0] k^-_i[0] = k^-_i[0] k^+_i[0] = 1,$$
$$k^\pm_i(z) k^\pm_j(w) = k^\pm_j(w) k^\pm_i(z),$$
$$k^-_i(z) k^+_i(w) = \frac{f(\frac{w}{z}q^{-k})}{f(\frac{w}{z}q^k)}
k^+_i(w) k^-_i(z),$$
$$k^-_1(z) k^+_2(w) = \frac{f(\frac{w}{z}q^{k+2})}{f(\frac{w}{z}q^{-k+2})}
k^+_2(w) k^-_1(z),$$
$$k^-_2(z) k^+_1(w) = \frac{f(\frac{w}{z}q^{k-2})}{f(\frac{w}{z}q^{-k-2})}
k^+_1(w) k^-_2(z),$$
$$k^\pm_1(z) E(w) = \frac{zq^{\mp\frac{k}{2}-1} - wq}{zq^{\mp\frac{k}{2}} -
w} E(w) k^\pm_1(z),$$
$$k^\pm_1(z) F(w) = \frac{zq^{\pm\frac{k}{2}} - w}{zq^{\pm\frac{k}{2}-1} -
wq} F(w) k^\pm_1(z),$$
$$k^\pm_2(z) E(w) = \frac{zq^{\mp\frac{k}{2}+1} -
wq^{-1}}{zq^{\mp\frac{k}{2}} - w} E(w) k^\pm_2(z),$$
$$k^\pm_2(z) F(w) = \frac{zq^{\pm\frac{k}{2}} - w}{zq^{\pm\frac{k}{2}+1} -
wq^{-1}} F(w) k^\pm_2(z),$$
$$E(z) E(w) = \frac{zq^2-w}{z-wq^2} E(w) E(z),$$
$$F(z) F(w) = \frac{z-wq^2}{zq^2-w} F(w) F(z),$$
\begin{multline*}
\left[ E(z),F(w) \right] = \\ (q-q^{-1}) \left( \delta \left(\frac{w}{z}
q^k \right) k^-_2(w q^{\frac{k}{2}}) k^-_1(w q^{\frac{k}{2}})^{-1} - \delta
\left(\frac{w}{z} q^{-k} \right) k^+_2(w q^{-\frac{k}{2}}) k^+_1(w
q^{-\frac{k}{2}})^{-1} \right),
\end{multline*}
where $$\delta(x) = \sum_{m\in\Z} x^m.$$

These relations are understood as relations between formal powers series
(cf. \remref{power}).

\begin{lem}
The Fourier coefficients of the power series $k^\pm_1(z)
k^\pm_2(zq^{-2})-1$ are central elements of $U_h(\gl)_k$.
\end{lem}

Consider the quotient of $U_h(\gl)_k$ by the ideal generated by these
elements. It has $E[n], F[n], k^\pm_1[n], n\in\Z$, as generators. There is
a one-to-one correspondence between them and Drinfeld's generators which
preserves relations: Drinfeld's $\psi(z)$ ($\psi_+(z)$ of \cite{AOS}) is
$k^-_1(zq^2)^{-1} k^-_1(z)^{-1}$, $\phi(z)$ ($\psi_-(z)$ of \cite{AOS}) is
$k^+_1(zq^2)^{-1} k^+_1(z)^{-1}$, $\xi^+(z)$ is $E(z)$ and $\xi^-(z)$ is
$h(q-q^{-1}) F(z)$. The following Proposition then follows from
\cite{Dr:new}, cf. also \cite{KhT,B,LSS}.

\begin{prop}    \label{gl2}
The quotient of $U_h(\gl)_k$ by the ideal generated by the Fourier
coefficients of the power series $k^\pm_1(z) k^\pm_2(zq^{-2})-1$ is
isomorphic to $U_h(\su)_k$.
\end{prop}

\section{The RS realization.}
Now we turn to the RS realization \cite{RS}. This realization originated
from the Quantum Inverse Scattering Method, cf. \cite{Fad,FRT}. It is
important for us because in this realization we can write explicit formulas
for central elements \cite{RS}.

Introduce the $R$--matrix
\begin{equation}    \label{Rmatrix}
R(x) = f(x) \begin{pmatrix}
1 & 0 & 0 & 0 \\
0 & \dfrac{1-x}{q-xq^{-1}} & \dfrac{x(q-q^{-1})}{q-xq^{-1}} & 0 \\
0 & \dfrac{q-q^{-1}}{q-xq^{-1}} & \dfrac{1-x}{q-xq^{-1}} & 0 \\
0 & 0 & 0 & 1
\end{pmatrix}
\end{equation}
where $f(x)$ is given by formula \eqref{f}.

The matrix \eqref{Rmatrix} is the result of computation of the universal
$R$--matrix of $U_q(\su)$ on the tensor product of two two-dimensional
evaluation representations, cf. e.g. \cite{FR}. It satisfies the
crossing-symmetry property:
$$\left( \left( \left( R(x)^{-1} \right)^{t_1} \right)^{-1} \right)^{t_1} =
\left( \begin{pmatrix}
q^{-1} & 0 \\
0 & q
\end{pmatrix} \otimes I_2 \right) R(xq^4) \left( \begin{pmatrix}
q & 0 \\
0 & q^{-1}
\end{pmatrix} \otimes I_2 \right),$$
which follows from the existence of an isomorphism between the
two-dimensional evaluation module and its double dual \cite{FR}. This is
related to the fact that $f(x)$ satisfies the $q$-difference equation
$$f(xq^4) = \frac{(1-xq^2)^2}{(1-x)(1-xq^4)} f(x).$$

\begin{rem}
Our $R$--matrix differs from that of \cite{DF} by the factor $f(x)$ and
by replacement of $q$ by $q^{-1}$. It also differs from the $R$--matrix
used in \cite{JM} by the factor which is a product of theta-functions.
\qed\end{rem}

Let again $h \in \C \backslash \{ 2\pi i {\Bbb Q} \}, k \in \C$, $q=e^h$.
We define an associative algebra $U'_h(\gl)_k$ over $\C$ with generators
$l^\pm_{ij}[n]$, where $i,j=1,2$, and $n\in \mp\Z_+\backslash 0$, and
$l^+_{ij}[0], l^-_{ji}[0], 1\leq j \leq i \leq 2$. Introduce the generating
functions $$l^\pm_{ij}(z) = \sum_{m=0}^\infty l^\pm_{ij}[\mp n] z^{\pm n},$$
where we put $l^+_{ij}[0]=l^-_{ji}[0]=0$ for $1\leq i < j \leq 2$. Let
$L^\pm(z)$ be the $2\times 2$ matrix $(l^\pm_{ij}(z))_{i,j=1,2}$.

The defining relations in $U'_h(\gl)_k$ are: $$l^+_{ii}[0] l^-_{ii}[0] =
l^-_{ii}[0] l^+_{ii}[0] = 1, \quad \quad i=1,2,$$
\begin{equation}    \label{first}
R \left( \frac{z}{w} \right) L^\pm_1(z) L^\pm_2(w) = L^\pm_2(w) L^\pm_1(z)
R \left( \frac{z}{w} \right),
\end{equation}
\begin{equation}    \label{second}
R \left( \frac{z}{w} q^{-k} \right) L^+_1(z) L^-_2(w) = L^-_2(w) L^+_1(z) R
\left( \frac{z}{w} q^k \right).
\end{equation}
Here $L^\pm_1(z)$ and $L^\pm_2(w)$ are elements of $\on{End}(\C^2)\otimes
\on{End}(\C^2) \otimes U'_h(\gl)_k$, i.e. $4\times 4$ matrices with entries
from $U'_h(\gl)_k$, such that $L^\pm_1(w)=L^\pm(w) \otimes I_2$,
$L^\pm_2(z)=I_2 \otimes (L^\pm(z))$.

The relations \eqref{first} and \eqref{second} are understood as relations
between formal power series in $\dfrac{z}{w}$, cf. \remref{power}.

\section{The isomorphism of two realizations.}
Following \cite{DF} one can construct an explicit isomorphism between the
algebras $U_h(\gl)_k$ and $U'_h(\gl)_k$.

Consider the following decomposition
$$L^\pm(z) =
\begin{pmatrix}
1 & 0 \\
e^\pm(z) & 1
\end{pmatrix}
\begin{pmatrix}
k^\pm_1(z) & 0 \\
0 & k^\pm_2(z)
\end{pmatrix}
\begin{pmatrix}
1 & f^\pm(z) \\
0 & 1
\end{pmatrix}$$
\begin{equation}    \label{Lop}
=
\begin{pmatrix}
k^\pm_1(z) & k^\pm_1(z) f^\pm(z) \\
e^\pm(z) k^\pm_1(z) & k^\pm_2(z) + e^\pm(z) k^\pm_1(z) f^\pm(z)
\end{pmatrix}.
\end{equation}

In particular, we see that $l^+_{21}[0]$ is the constant term of
$e^+(z)=\sum_{m\geq 0} e[-m] z^m$ and $l^-_{12}[0]$ is the constant term of
$f^-(z) = \sum_{m\leq 0} f[-m] z^m$, while $e^-(z) = \sum_{m<0} e[-m] z^m$
and $f^+(z) = \sum_{m>0} f[-m] z^m$ have no constant terms.

\begin{prop}[\cite{DF}] The map $\psi': U'_h(\gl)_k \arr U_h(\gl)_k$
defined on generators by
$$\psi'[k^\pm_1(z)] = k^\pm_1(z),$$ $$\psi'[k^\pm_2(z)]=
k^\pm_2(z),$$ $$\psi'[e^+(zq^{\frac{k}{2}}) -
e^-(zq^{-\frac{k}{2}})] = E(z),$$ $$\psi'[f^+(zq^{-\frac{k}{2}}) -
f^-(zq^{\frac{k}{2}})] = F(z),$$ is an isomorphism.
\end{prop}

\begin{lem}
The Fourier coefficients of the power series
\begin{equation}    \label{diagonal}
l^\pm_{11}(zq^2) \left( l^\pm_{22}(z) - l^\pm_{21}(z) l^\pm_{11}(z)^{-1}
l^\pm_{12}(z) \right) -1
\end{equation}
are central elements of $U'_h(\gl)_k$.
\end{lem}

The image of the power series \eqref{diagonal} under the map $\psi$ is the
series $k^\pm_1(zq^2) k^\pm_2(z)-1$, whose coefficients are central
elements of $U_h(\gl)_k$. According to \propref{gl2}, the quotient of
$U_h(\gl)_k$ by the ideal generated by these central elements is isomorphic
to $U_h(\su)_k$. Hence we obtain

\begin{cor} $\psi'$ induces an isomorphism $\psi$ between the quotient of
$U'_h(\gl)_k$ by the Fourier coefficients of the power series
$l^\pm_{11}(zq^2) ( l^\pm_{22}(z) - l^\pm_{21}(z) l^\pm_{11}(z)^{-1}
l^\pm_{12}(z) ) - 1$ and $U_h(\su)_k$.
\end{cor}

\section{The $q$--analogues of the Sugawara operators.}
We define a completion $\wt{U}_h(\G)_k$ of $U_h(\G)_k$ as follows:
$$\wt{U}_h(\G)_k = \underset{\longleftarrow}{\lim} \; U_h(\G)_k/J_n, \quad
\quad n>0,$$ where $J_n$ is the left ideal of $U_h(\G)_k$ generated by
$l^-_{ij}[m], m\geq n$.

Let $$L(z) = L^+(zq^{-\frac{k}{2}}) L^-(zq^{\frac{k}{2}})^{-1}.$$ It is
easy to see that all Fourier coefficients of the power series
\begin{equation}    \label{central}
\el(z) = q^{-1} L_{11}(z) + q L_{22}(z),
\end{equation}
lie in $\wt{U}_h(\su)_k$. It follows from \cite{RS} that when $k=-2$ the
coefficients of $\el(z)$ are central elements of $\wt{U}_h(\su)\cri$.

We will now express $\el(z)$ in terms of the Drinfeld realization using the
isomorphism $\psi$. Let us put $k=-2$. Using formula \eqref{Lop} we obtain:
\begin{equation}    \label{L11}
L_{11}(z) = k^+_1(zq) k^-_1(zq^{-1})^{-1} - k^+_1(zq) \left( f^+(zq) -
f^-(zq^{-1}) \right) k^-_2(zq^{-1})^{-1} e^-(zq^{-1}),
\end{equation}
\begin{equation}    \label{L22}
L_{22}(z) = k^+_2(zq) k^-_2(zq^{-1})^{-1} + e^+(zq) k^+_1(zq) \left(
f^+(zq) - f^-(zq^{-1}) \right) k^-_2(zq^{-1})^{-1}.
\end{equation}

Applying formula (4.45) from \cite{DF} (in which the sign of the second
summand in the right hand side has to be reversed) we obtain
\begin{equation}    \label{DF1}
k^-_2(zq^{-1})^{-1} e^-(zq^{-1}) = q e^-(zq) k^-_2(zq^{-1})^{-1},
\end{equation}
and applying formula (4.25) from \cite{DF} we obtain
\begin{equation}    \label{DF2}
e^+(zq) k^+_1(zq) = q^{-1} k^+_1(zq) e^+(zq^{-1}).
\end{equation}
Substituting \eqref{DF1} into \eqref{L11} we obtain
\begin{equation}    \label{mL11}
L_{11}(z) = k^+_1(zq) k^-_1(zq^{-1})^{-1} - q k^+_1(zq) \left( f^+(zq) -
f^-(zq^{-1}) \right) e^-(zq) k^-_2(zq^{-1})^{-1},
\end{equation}
and substituting \eqref{DF2} into \eqref{L22} we obtain
\begin{equation}    \label{mL22}
L_{22}(z) = k^+_2(zq) k^-_2(zq^{-1})^{-1} + q^{-1} k^+_1(zq) e^+(zq^{-1})
\left( f^+(zq) - f^-(zq^{-1}) \right) k^-_2(zq^{-1})^{-1}.
\end{equation}

Combining formula \eqref{central} with \eqref{mL11} and \eqref{mL22} we
obtain
\begin{align}    \label{t1}
\el(z) &= q^{-1} k^+_1(zq) k^-_1(zq^{-1})^{-1} + q k^+_2(zq)
k^-_2(zq^{-1})^{-1} \\    \label{t2} &+ k^+_1(zq) :E(z) F(z):
k^-_2(zq^{-1})^{-1},
\end{align}
where $$:E(z) F(z): = e^+(zq^{-1}) \left( f^+(zq) - f^-(zq^{-1}) \right)
- \left( f^+(zq) - f^-(zq^{-1}) \right) e^-(zq).$$

Now we apply the isomorphism $\psi$. We see that $\psi(:E(z) F(z):)$ is a
normally ordered product of the power series $E(z)$ and $F(z)$: $$:E(z)
F(z): = E_-(z) F(z) + F(z) E_+(z),$$ where $$E_-(z) = \sum_{n\leq 0} E[n]
z^{-n}, \quad \quad E_+(z) = \sum_{n\geq 0} E[n] z^{-n}.$$ Thus, we obtain

\begin{prop} The Fourier coefficients of the power series
\begin{align}    \label{qsug}
\el(z) &= q^{-1} k^+_1(zq) k^-_1(zq^{-1})^{-1} + q k^+_1(zq^3)^{-1}
k^-_1(zq) \\ &+ k^+_1(zq) :E(z) F(z): k^-_1(zq)
\end{align}
are central elements of $\wt{U}_h(\su)\cri$.
\end{prop}
These are the $q$--analogues of the Sugawara operators.

\section{The Wakimoto realization of $U_h(\su)_k$.}
Now we will describe a homomorphism $\phi_{h,k}$ from $U_h(\su)_k$ to a
completion of a quantum Heisenberg algebra. The map $\phi_{h,k}$ is a
$q$--analogue of the map $\phi$ defined in \S~2. Such homomorphisms have
been constructed in \cite{qWak}. We will use the Awata--Odake--Shiraishi
construction \cite{AOS} with some modifications.

Introduce the quantum Heisenberg algebra ${\cal A}_{h,k}(\sw_2)$. The
generators are $\la_n, b_n, c_n$, $n\in\Z, n\neq 0$, $\exp(\pm \la_0/2),
\exp(\pm (q-q^{-1})b_0/2), \exp(\pm (q-q^{-1})c_0)$, and $p_b$, $p_c$. The
relations are
\begin{equation}    \label{lrel}
[\la_n,\la_m] = \frac{1}{n} \frac{[(k+2)n]_q [n]_q^2}{[2n]_q}
\delta_{n,-m},
\end{equation}
$$[b_n,b_m] = -\frac{1}{n} [n]_q^2 \delta_{n,-m}, \quad \quad [b_0,p_b] = -
\frac{q-q^{-1}}{2 h}$$
$$[c_n,c_m] = \frac{1}{n} [n]_q^2 \delta_{n,-m}, \quad \quad [c_0,p_c] =
\frac{q-q^{-1}}{2 h}$$ where $q=e^h$ and $[n]_q =
\dfrac{q^n-q^{-n}}{q-q^{-1}}$. The generators $b_n$ and $c_n$ coincide with
the corresponding generators of \cite{AOS}. The generators $\la_n$ are
related to the generators $a_n$ of \cite{AOS} by the formula $$\la_n =
\frac{q-q^{-1}}{q^n+q^{-n}} a_n$$ (recall that we have assumed that $q$ is
not a root of unity).

We form the generating functions:
$$\La_\pm(z) = \exp \left( \pm \frac{\la_0}{2} \pm \sum_{n=1}^\infty
\la_{\pm n} z^{\mp n} \right),$$
$$b_\pm(z) = \pm (q-q^{-1}) \left( \frac{b_0}{2} + \sum_{n=1}^\infty
b_{\pm n} z^{\mp n} \right),$$
$$b(z) = - \sum_{n\neq 0} \frac{b_n}{[n]_q} z^{-n} + \frac{q-q^{-1}}{2h}
b_0 \log z + p_b,$$
$$c_\pm(z) = \pm (q-q^{-1}) \left( \frac{c_0}{2} + \sum_{n=0}^\infty c_{\pm
n} z^{\mp n}, \right),$$
$$c(z) = - \sum_{n\neq 0} \frac{c_n}{[n]_q} z^{-n} + \frac{q-q^{-1}}{2h}
c_0 \log z + p_c.$$

The series $\La_\pm(z)$ is related to the series $a_{\pm}(z)$ from
\cite{AOS} by the formula
$$\La_\pm(z) \La_\pm(zq^{\pm 2}) = e^{a_\pm(z q^{\pm 1})}.$$ The other
series are the same as in \cite{AOS}.

The relations between these series, in the sense of formal power series
(cf. \remref{power}), are the following (cf. \cite{AOS}):
$$\La_+(z) \La_-(w) = \frac{f(\frac{w}{z} q^{-k-2})}{f(\frac{w}{z}
q^{k+2})} \La_-(w) \La_+(z),$$
\begin{equation}    \label{r1}
e^{b_+(z)} :e^{b(w)}: = \frac{z-wq}{zq-w} :e^{b(w)}: e^{b_+(z)},
\end{equation}
\begin{equation}    \label{r2}
:e^{b(z)}: e^{b_-(w)} = \frac{z-wq}{zq-w} e^{b_-(w)} :e^{b(z)}: = q
\frac{z-wq}{zq-w} :e^{b_-(w)+b(z)}:.
\end{equation}
\begin{equation}    \label{r3}
e^{b_+(z)} e^{b_-(w)} = \frac{(z-w)^2}{(z-wq^2)(z-wq^{-2})} e^{b_-(w)}
e^{b_+(z)},
\end{equation}
$$e^{c_+(z)} e^{c_-(w)} = \frac{(z-wq^2)(z-wq^{-2})}{(z-w)^2} e^{c_-(w)}
e^{c_+(z)},$$

We define the completion $\ca$ of ${\cal A}_{h,k}(\sw_2)$ as follows:
$$\ca = \underset{\longleftarrow}{\lim} \; {\cal A}_{h,k}/I_n, \quad \quad
n>0,$$ where $I_n$ is the left ideal of ${\cal A}_{h,k}(\sw_2)$ generated
by all polynomials in $\la_m,b_m,c_m, m>0$, of degrees greater than or
equal to $n$ (we put $\deg \la_m = \deg b_m = \deg c_m = m$).

The next proposition follows from \cite{AOS}.

\begin{prop}
There is a homomorphism $\phi_{h,k}$ from $U_h(\su)_k$ to $\ca$, which is
defined on generators as follows:
\begin{align*}
\phi_{h,k}[E(z)] &= -:e^{b_+(z)-(b+c)(zq)}: + :e^{b_-(z)-(b+c)(zq^{-1})}:, \\
\phi_{h,k}[F(z)] &= \La_+(zq^{\frac{k}{2}}) \La_+(zq^{\frac{k}{2}+2})
:e^{b_+(zq^{k+2})+(b+c)(zq^{k+1})}:, \\ &- \La_-(zq^{-\frac{k}{2}})
\La_-(zq^{-\frac{k}{2}-2})  :e^{b_-(zq^{-k-2})+(b+c)(zq^{-k-1})}:,\\
\phi_{h,k}[k^+_1(z)] &= \La_-(zq^{-2})^{-1} e^{-b_-(zq^{-\frac{k}{2}-2})}, \\
\phi_{h,k}[k^-_1(z)] &= \La_+(z)^{-1} e^{-b_+(zq^{\frac{k}{2}})}.
\end{align*}
\end{prop}

\begin{rem} Under the homomorphism $\phi$ defined in \S~2.2, the affine
algebra $\su$ embeds into a completion of the Heisenberg algebra $\Gamma
\otimes {\cal H}'(\sw_2)$ generated by $a_n, a_n^*$ and $\chi_n$. The power
series $a(z)$ and $a^*(z)$ form the so-called $\beta\gamma$--system while
the power series $\chi(z)$ is called a free scalar field. The
$\beta\gamma$--system can be represented via exponentials of a pair of free
scalar fields. The homomorphism $\phi$ then gives rise to a homomorphism
$\phi'$ from $\su$ to a completion of the Heisenberg algebra generated by
the Fourier coefficients of these two scalar fields and $\chi(z)$.

The power series $b(z)$ and $c(z)$ are $q$--analogues of the scalar bosonic
fields representing the $\beta\gamma$ system when $q=1$. Thus, the
homomorphism $\phi_{h,k}$ is a $q$--deformation of $\phi'$ rather than
$\phi$.
\qed\end{rem}

When $k\neq -2$, the homomorphism $\phi_{h,k}$ provides representations of
$U_h(\su)_k$ in the Fock representation of the Heisenberg algebra ${\cal
A}_{h,k}$, cf. \cite{AOS}. These representations have one parameter -- the
action of $\la_0$ on the highest weight vector, cf. \cite{AOS}. When
$k=-2$, the generators $\la_n$ commute among themselves and generate a
commutative algebra ${\cal H}_h(\sw_2)$. Therefore representations of
$U_h(\su)\cri$ can be realized via $\phi_{h,k}$ in a smaller space: the
tensor product of the Fock representation of the subalgebra of ${\cal
A}_{h,-2}$ generated by $b_n, c_n, n\in\Z$, and a one-dimensional
representation of ${\cal H}_h(\sw_2)$. For the action of
$\wt{U}_h(\su)\cri$ to be well-defined on this space, the action of
$\La(z)$ should be given by a Laurent power series $\la(z)$ (for example,
this is the case if $\la_n, n>0$, act by $0$). The corresponding
representations $W_{\la(z)}$ are the $q$--analogues of the Wakimoto modules
over $U(\su)\cri$ from \S~2.2.

\section{Deformation of the Miura transformation.}
In this section we will apply $\phi_h^{\on{cr}} \equiv \phi_{h,-2}$ to the
generating function of central elements $\el(z)$ given by \eqref{qsug}. For
brevity, in what follows we will use the same notation for elements of
$\wt{U}_h(\su)\cri$ and their images in $\ca$.

The normally ordered product $:E(z) F(z):$ can be written as
\begin{equation}    \label{normal}
:E(z) F(z): = \int_{C_R} \frac{E(w) F(z)}{w-z} dw - \int_{C_r} \frac{F(z)
E(w)}{w-z},
\end{equation}
where $C_R$ and $C_r$ are circles around the origin or radii
$R>|w|$ and $r<|w|$, respectively.

Using the Wakimoto realization of
$U_h(\su)\cri$ and formulas \eqref{r1}--\eqref{r3}, we find in the region
$|w|>|z|, |w|>q^2|z|, |w|>q^{-2}|z|$:
\begin{align*}
E(w) F(z) = &- \frac{w-z}{wq-zq^{-1}} \La_+(zq^{-1})
\La_+(zq) e^{(b+c)(zq^{-1}) - (b+c)(wq)} e^{b_+(z) + b_+(w)} \\
&- \frac{w-z}{wq^{-1}-zq} \La_-(zq^{-1}) \La_-(zq) e^{b_-(z) +
b_-(w)} e^{(b+c)(zq)-(b+c)(wq^{-1})} \\
&+ q^{-1} \La_-(zq^{-1}) \La_-(zq) e^{b_-(z)} e^{(b+c)(zq) -
(b+c)(wq)} e^{b_+(w)}  \\
&+ q \La_+(zq^{-1}) \La_+(zq) e^{b_-(w)}
e^{(b+c)(zq^{-1}) - (b+c)(wq^{-1})} e^{b_+(z)}.
\end{align*}

We obtain the same formula for $F(z) E(w)$ in the region $|w|<|z|,
|w|<q^2|z|, |w|<q^{-2}|z|$. We can therefore rewrite \eqref{normal} as the
integral of the this expression over the contour on the $w$ plane
surrounding the points $z, zq^2, zq^{-2}$.

Evaluating the residues, we find that
\begin{align}    \label{norm1}
:E(z) F(z): = &- q^{-1} \La_+(zq^{-1}) \La_+(zq) e^{b_+(zq^{-2}) + b_+(z)}
\\    \label{norm2} &-
q \La_-(zq^{-1}) \La_-(zq) e^{b_-(zq^2)+b_-(z)} \\    \label{norm3}
&+ q^{-1} \La_-(zq^{-1}) \La_-(zq)e^{b_-(z)} e^{b_+(z)} \\    \label{norm4}
&+ q \La_+(zq^{-1}) \La_+(zq) e^{b_-(z)} e^{b_+(z)}.
\end{align}

Using this formula we obtain
\begin{align*}
k^+_1(zq) :E(z) F(z): k^-_1(zq) = & - q^{-1} \La_-(zq^{-1})^{-1}
\La_+(zq^{-1}) e^{-b_-(z)} e^{b_+(zq^{-2})} \\ & - q \La_-(zq)
\La_+(zq)^{-1} e^{b_-(zq^2)} e^{-b_+(z)} \\ & + q^{-1} \La_-(zq)
\La_+(zq)^{-1} \\ & + q \La_-(zq^{-1})^{-1} \La_+(zq^{-1}).
\end{align*}

On the other hand, we have $$k^+_1(zq) k^-_1(zq^{-1})^{-1} =
\La_-(zq^{-1})^{-1} \La_+(zq^{-1}) e^{-b_-(z)} e^{b_+(zq^{-2})}$$ and
$$k^+_1(zq^3)^{-1} k^-_1(zq) = \La_-(zq) \La_+(zq)^{-1} e^{b_-(zq^2)}
e^{-b_+(z)}.$$ Substituting these formulas into \eqref{qsug} we obtain a
formula expressing the image of $\el(z)$ in ${\cal H}_h(\sw_2)$ in terms
of $\La_\pm(z)$.

\begin{thm} {\em Under the homomorphism $\phi_h^{\on{cr}}$,}
\begin{equation}    \label{final}
\el(z) \longrightarrow \s(z) = \La(zq) + \La(zq^{-1})^{-1},
\end{equation}
{\em where}
\begin{equation}    \label{Lambda}
\La(z) = q^{-1} \La_-(z) \La_+(z)^{-1}.
\end{equation}
\end{thm}

This is a $q$--deformation of the Miura transformation \eqref{mm}.

\begin{rem}    \label{simple}
There is a simpler way to obtain formula \eqref{final}. Consider the action
of $\s(z) = \phi_h^{\on{cr}}[\el(z)]$ on the module $W_{\la(z)}$ introduced
at the end of \S~7. In the limit $q \arr 1$ this module becomes a Wakimoto
module over $\su$. Wakimoto modules are irreducible for generic values of
parameters. Therefore the same is true for the modules $W_{\la(z)}$. Hence
any central element of $\wt{U}_h(\su)\cri$ acts on $W_{\la(z)}$ by a
constant. In particular, $\el(z)$ acts by multiplication by a Laurent power
series $\wt{\s}(z)$. We can compute $\wt{\s}(z)$ by taking the matrix
element of $\el(z)$ between the generating vector of $W_{\la(z)}$ and its
dual using formulas \eqref{central}, \eqref{Lop} and the maps $\psi,
\phi_h^{\on{cr}}$. Explicit computation shows that this matrix element is
equal to $\la(zq)+\la(zq^{-1})^{-1}$. This implies that $\s(z)$ lies in
${\cal H}_h(\sw_2)$ and gives us formula \eqref{final} for $\s(z)$.
\qed\end{rem}

\begin{rem}
In \cite{RS} a generating function of central elements of
$\wt{U}_h(\G)\cri$ has been associated to an arbitrary finite-dimensional
representation of $U_q(\G)$. Thus, in the case of $\wt{U}_h(\su)\cri$ we
have a generating function $\el^{(n)}(z)$ of central elements associated to
the representation of $U_q(\su)$ of dimension $n+1$ for each positive
integer $n$. In particular, $\el(z) = \el^{(1)}(z)$. These generating
functions satisfy the following relation:
$$\el^{(1)}(zq^{2n}) \el^{(n)}(z) = \el^{(n+1)}(z) + \el^{(n-1)}(z), \quad
\quad n>0,$$ Using this relation and formula \eqref{final} it is easy to
find recursively:
\begin{align*}
\el^{(n)}(z) &= \La(zq) \La(zq^3) \La(zq^5) \ldots \La(zq^{2n-1}) \\
             &+ \La(zq^{-1})^{-1} \La(zq^3) \La(zq^5) \ldots \La(zq^{2n-1})
\\           &+ \La(zq^{-1})^{-1} \La(zq)^{-1} \La(zq^5) \ldots
\La(zq^{2n-1}) \\
             &+ \La(zq^{-1})^{-1} \La(zq)^{-1} \La(zq^3)^{-1} \ldots
\La(zq^{2n-1}) \\
             &+ \ldots \\
             &+ \La(zq^{-1})^{-1} \La(zq)^{-1} \La(zq^3)^{-1} \ldots
\La(zq^{2n-3})^{-1}
\end{align*}
(compare with \cite{KR})
\qed\end{rem}

\section{Poisson bracket.}
We can now compute the Poisson bracket between Fourier coefficients of
$\s(z)$. They generate a central subalgebra of $\wt{U}_h(\su)\cri$. Let
$Z_h(\su)$ be its completion in $\wt{U}_h(\su)\cri$.

Consider the algebra $\Hh_{h,k}(\sw_2)$ with generators $\la_n, n\in\Z$,
$\exp(\pm \la_0)$ and relations \eqref{lrel}. Let $\wt{{\cal
H}}_{h,k}(\sw_2)$ be its completion defined as follows:
$$\wt{{\cal H}}_{h,k}(\sw_2) = \underset{\longleftarrow}{\lim} \; {\cal
H}_{h,k}(\sw_2)/I_n, \quad \quad n>0,$$ where $I_n$ is the left ideal of
$\Hh_{h,k}(\sw_2)$ generated by all polynomials in $\la_m, m>0$, of degrees
greater than or equal to $n$.

The family $\wt{{\cal H}}_{h,k}(\sw_2)$ induces a Poisson structure on ${\cal
H}_h(\sw_2) \equiv \wt{{\cal H}}_{k,-2}(\sw_2)$, such that
\begin{equation}    \label{skobka}
\{ \la_n,\la_m \} = 2h(q-q^{-1}) \frac{[n]_q^2}{[2n]_q} \delta_{n,-m}
\end{equation}
(recall that $q=e^h$ and that $h\not{\hspace*{-1mm}\in} \; 2\pi i {\Bbb Q}$).

We define a Poisson bracket on the center of $\wt{U}_h(\su)\cri$ in the
same way as in \S~2.1, as the leading term in the commutator of liftings of
central elements. Formula \eqref{final} shows that the images of the
Fourier coefficients of $\el(z)$ under the homomorphism $\phi_h^{\on{cr}}$
lie in ${\cal H}_h(\sw_2)$. As was explained in \S~2.2, we can take
liftings of the images of central elements inside the deformation
$\wt{{\cal H}}_{h,k}(\sw_2)$ of ${\cal H}_h(\sw_2)$. This shows that the
homomorphism from the center to ${\cal H}_h(\sw_2)$ is a homomorphism of
Poisson algebras. We can use this fact to compute the Poisson brackets
between Fourier coefficients of $\s(z)$.

{}From formulas \eqref{skobka} and \eqref{Lambda} we obtain
\begin{equation}    \label{lpb}
\{ \La(z),\La(w) \} = 2h (q-q^{-1}) \sum_{m\in\Z} \left( \frac{w}{z}
\right)^m \qn \La(z) \La(w).
\end{equation}

According to formula \eqref{final}, we have:
\begin{align*}
\{ \s(z),\s(w) \} &= \{ \La(zq),\La(wq) \} + \{
\La(zq),\La(wq^{-1})^{-1} \} \\ &+ \{
\La(zq^{-1})^{-1},\La(wq) \} +
\{ \La(zq^{-1})^{-1},\La(wq^{-1})^{-1} \}.
\end{align*}
Substituting \eqref{lpb} into this formula, we obtain:
\begin{align*}
\{ \s(z),\s(w) \} &= 2h(q-q^{-1}) \sum_{m\in\Z} \left( \frac{w}{z} \right)^m
\qn \La(zq) \La(wq) \\ &- 2h(q-q^{-1}) \sum_{m\in\Z} \left(
\frac{w}{z} \right)^m \qn q^{-2m} \La(zq)
\La(wq^{-1})^{-1} \\
&- 2h(q-q^{-1}) \sum_{m\in\Z} \left( \frac{w}{z} \right)^m
\qn q^{2m} \La(zq^{-1})^{-1} \La(wq) \\ &+ 2h(q-q^{-1})
\sum_{m\in\Z} \left( \frac{w}{z} \right)^m \qn
\La(zq^{-1})^{-1} \La(wq^{-1})^{-1} \\
&= 2h(q-q^{-1}) \sum_{m\in\Z} \left( \frac{w}{z} \right)^m \qn
\s(z) \s(w) \\
&- 2h \sum_{m\in\Z} \left( \frac{w}{z} \right)^m (1-q^{-2m})
\La(zq) \La(wq^{-1})^{-1} \\ &- 2h \sum_{m\in\Z} \left( \frac{w}{z}
\right)^m (q^{2m}-1) \La(zq^{-1})^{-1} \La(wq).
\end{align*}
The last two terms give us $$- 2h \delta \left( \frac{w}{z} \right) \La(zq)
\La(wq^{-1})^{-1} + 2h \delta \left( \frac{w}{zq^2} \right) \La(zq)
\La(wq^{-1})^{-1}$$ $$- 2h \delta \left( \frac{wq^2}{z} \right)
\La(zq^{-1})^{-1} \La(wq) + 2h \delta \left( \frac{w}{z} \right)
\La(zq^{-1})^{-1} \La(wq)$$ $$= 2h \delta \left( \frac{w}{zq^2} \right) -
2h \delta \left( \frac{wq^2}{z} \right),$$ where $$\delta(x) =
\sum_{m\in\Z} x^m.$$

Finally, we obtain

\begin{thm}
\begin{multline}    \label{finalpb}
\{ \s(z),\s(w) \} = \\ 2h (q-q^{-1}) \sum_{m\in\Z} \left( \frac{w}{z}
\right)^m \qn \; \s(z) \s(w) +  2h \delta \left( \frac{w}{zq^2}
\right) - 2h \delta \left( \frac{wq^2}{z} \right).
\end{multline}
\end{thm}

This implies that the Poisson bracket between $\el(z)$ and $\el(w)$ is
given by formula \eqref{mult}.

Formula \eqref{finalpb} gives us the following formula for the Poisson
bracket between the Fourier coefficients of $\s(z)$:
\begin{equation}    \label{comp}
\{ \s_n,\s_m \} = 2h \sum_{l\in\Z} \frac{q^l-q^{-l}}{q^l+q^{-l}}
\s_{n-l} \s_{m+l} - 2h (q^{2n}-q^{-2n}) \delta_{n,-m},
\end{equation}
where we put $$\s(z) = \sum_{n\in\Z} \s_n z^{-n}.$$

The elements $\s_n$ generate a Poisson algebra $\W_h(\sw_2)$, which is a
$q$--deformation of the classical Virasoro algebra $\W(\sw_2)$. This
Poisson algebra is embedded into ${\cal H}_h(\sw_2)$ via the
$q$--deformation of the Miura transformation \eqref{final}. The Poisson
algebra $Z_h(\su)$ is isomorphic to $\W_h(\sw_2)$.

\begin{rem}
In the limit $h \arr 0$ we have: $$\el(z) = 2 + 4h^2\left( z^2 S(z) +
\frac{1}{4} \right) + h^3(\ldots)$$ (cf. \cite{DE}). On the other hand
$$\La(z) = 1 - h (z\chi(z)+1) + h^2(\ldots).$$

Substituting these formulas into formula \eqref{final} and expanding in
powers of $h$ up to $h^2$, we obtain:
$$2 + 4h^2 z^2 S(z) + h^2 = 2 + h^2 \left( z\chi(z) + 1 \right)^2 - 2 h^2
z \pa_z (z\chi(z)),$$ which coincides with the Miura transformation
\eqref{mm}.

Now let us consider formula \eqref{comp}. The leading term in the expansion
of the left hand side is $16 h^4 \{ S_n,S_m \}$. Expanding the right hand
side of \eqref{comp} in powers of $h$ up to $h^4$, we obtain for the first
term $$8 h^2 n \delta_{n,-m} + 16 h^4 (n-m) \left( S_{n+m}  + \frac{1}{4}
\delta_{n,-m} \right) - h^4 \frac{8}{3} n^3 \delta_{n,-m},$$ and for the
second term $$- 8 h^2 n \delta_{n,-m} - h^4 \frac{16}{3} n^3
\delta_{n,-m}.$$ Taking the sum, we see that the leading term in the right
hand side is $16h^4$ times the right hand side of formula \eqref{sugpb}.

Note that we can obtain a different Poisson algebra by placing an arbitrary
overall factor in the right hand side of formula \eqref{finalpb}. In
particular, if we put the overall factor $-c/6$ in the right hand side of
the formula, then in the limit $h \arr 0$ we will recover the classical
Virasoro algebra with central charge $c$.

We can also replace the overall factor $2h$ in the right hand side by
$(q-q^{-1})$ without changing the asymptotics $h \arr 0$. After that we can
consider $\W_h(\sw_2)$ as a Poisson algebra over the ring of rational
functions in $q$.
\qed\end{rem}

\begin{rem}    \label{dual}
Formula \eqref{finalpb} gives another asymptotics as $h \arr 0$, if we
postulate that $\s(z)$ does not depend on $h$. If we divide the right hand
side of the formula by $2h(q-q^{-1})$, we obtain for $h=0$:
$$\{ \s(z),\s(w) \} = \frac{1}{2} \delta' \left( \frac{w}{z} \right) \s(z)
\s(w) - 2 \delta' \left( \frac{w}{z} \right).$$ The corresponding limiting
Poisson algebra $\wt{\W}(\sw_2)$ has a nice interpretation in terms of the
algebra ${\cal A}_k(\widehat{SL}_2^*)$ of functions on the Poisson-Lie
group $\widehat{SL}_2^*$ dual to $\widehat{SL}_2$ (this algebra is the
classical limit of $U_h(\su)_k$, cf. \cite{RS}). Namely, the limits of the
coefficients of $\el(z)$ generate a central Poisson subalgebra of ${\cal
A}_0(\widetilde{SL}_2^*)$ \cite{RS}. Deforming the level $k$ we obtain a
new Poisson structure on this subalgebra, which coincides with
$\wt{\W}(\sw_2)$. We will discuss this in more detail in \cite{next}.
\qed\end{rem}

\section{Generalization to $U_q(\sun)$.}

\subsection{A $q$--deformation of the Heisenberg-Poisson algebra.}
For a simple Lie algebra $\g$ of rank $l$, denote by
$B=(B_{ij})_{i,j=1,\ldots,l}$ the symmetrized Cartan matrix of $\g$; recall
that $B_{ij} = (\al_i,\al_j)$. Let ${\cal H}_{h,k}(\g)$ be the Heisenberg
algebra with generators $a_i[n], i=1,\ldots,l; n\in\Z$, and relations
$$[a_i[n],a_i[m]] = \frac{1}{n} [(k+\k)n]_q [B_{ij} n]_q \delta_{n,-m}.$$
The algebra ${\cal H}_{h,k}(\g)$ appears in \cite{AOS} in the construction
of the Wakimoto realization of $U_h(\sun)_k$.

The family ${\cal H}_{h,k}(\g)$ induces a Poisson structure on the
commutative algebra \newline ${\cal H}_{h,-\k}(\g)$. The Poisson brackets
between the generators of ${\cal H}_{h,-\k}(\g)$ are
\begin{equation}    \label{pba}
\{ a_i[n],a_j[m] \} = \frac{2h}{q-q^{-1}} [B_{ij} n]_q \delta_{n,-m}.
\end{equation}

Let ${\cal H}_h(\g)$ be the completion of ${\cal H}_{h,-\k}(\g)$ defined
in the same way as ${\cal H}_h(\sw_2)$.

Consider the case $\g=\sw_N$. Introduce new generators $\la_i[n],
i=1,\ldots,N; n\in\Z$, of ${\cal H}_h(\sw_N)$, which are related to the
generators $a_i[n]$ by the formulas
\begin{equation}    \label{conn}
\la_i[n] - \la_{i+1}[n] = q^{ni} (q-q^{-1}) a_i[n], \quad \quad
i=1,\ldots,N-1; n\in\Z,
\end{equation}
and which satisfy the linear relation
\begin{equation}    \label{linrel}
\sum_{i=1}^N q^{2(1-i)n} \la_i[n] = 0.
\end{equation}

We find from formulas \eqref{conn} and \eqref{linrel} the inverse change of
variables: $$\la_1[n] = (q-q^{-1}) \sum_{j=1}^{N-1}
\frac{[(N-j)n]_q}{[Nn]_q}
a_j[n],$$ $$\la_2[n] = -q^{Nn}(q-q^{-1}) \frac{[n]_q}{[Nn]_q} a_1[n] +
(q-q^{-1}) \sum_{j=2}^{N-1} \frac{[(N-j)n]_q}{[Nn]_q} a_j[n],$$ $$\ldots$$
$$\la_N[n] = - q^{Nn} (q-q^{-1}) \sum_{j=1}^{N-1} \frac{[jn]_q}{[Nn]_q}
a_j[n].$$ From these formulas and the brackets \eqref{pba} we find
\begin{align}    \label{pbl1}
\{ \la_i[n],\la_i[m] \} &= 2h (q-q^{-1}) \frac{[(N-1)n]_q [n]_q}{[Nn]_q}
\delta_{n,-m}, \\    \label{pbl2}
\{ \la_i[n],\la_j[m] \} &= - 2h (q-q^{-1}) \frac{[n]_q^2}{[Nn]_q} q^{-Nn}
\delta_{n,-m}, \quad \quad i<j.
\end{align}

Introduce the generating functions
\begin{equation}    \label{lai}
\La_i(z) = q^{-N+2i-1} \exp \left( - \sum_{m\in\Z} \la_i[m] z^{-m} \right).
\end{equation}
{}From \eqref{pbl1} and \eqref{pbl2}
we find:
\begin{align}    \label{pbg1}
\{ \La_i(z),\La_i(w) \} &= 2h (q-q^{-1}) \sum_{m\in\Z} \left( \frac{w}{z}
\right)^m \frac{[(N-1)m]_q [m]_q}{[Nm]_q} \La_i(z) \La_i(w),
\\    \label{pbg2}
\{ \La_i(z),\La_j(z) \} &= - 2h (q-q^{-1}) \sum_{m\in\Z} \left(
\frac{w}{zq^N} \right)^m \frac{[m]_q^2}{[Nm]_q} \La_i(z) \La_j(z), \quad
i<j.
\end{align}

\subsection{A $q$--deformation of the $\W$--algebra.}
Let us define generating functions $\s_i(z), i=0,\ldots,N$, whose
coefficients lie in $\HH$: $\s_0 = 1$, and
\begin{equation}    \label{qmap}
\s_i(z) = \sum_{1\leq j_1 < \ldots < j_i\leq N} \La_{j_1}(z) \La_{j_2}(zq^2)
\ldots \La_{j_{i-1}}(zq^{2(i-2)}) \La_{j_i}(zq^{2(i-1)}),
\end{equation}
$i=1,\ldots,N$. In particular,
$$\s_1(z) = \sum_{j=1}^N \La_j(z),$$ $$\s_2(z) = \sum_{1\leq j_1 < j_2 \leq
N} \La_{j_1}(z) \La_{j_2}(zq^2),$$ etc. Formula \eqref{linrel} implies that
$$\s_N(z) = \La_1(z) \La_2(zq^2) \ldots \La_N(zq^{2N-2}) = 1.$$

These formulas coincide with formulas for spectra of transfer-matrices in
integrable models associated to $U_q(\sun)$ \cite{KuR,BR}. Note that for
$\su$ we have $\s_1(z)=\s(zq^{-1})$, where $\s(z)$ is given by formula
\eqref{final}.

The coefficients of the series $\s_i(z), i=1,\ldots,N-1$, generate a
Poisson subalgebra $\W_h(\sw_N)$ of $\HH$. The relations between them can
be computed directly from formulas \eqref{pbg1} and \eqref{pbg2}. Introduce
the functions
\begin{equation}    \label{cij}
C_{ij}(x) = \sum_{m\in\Z} C^{(m)}_{ij} x^m, \quad \quad
i,j=1,\ldots,N-1,
\end{equation}
where
\begin{equation}    \label{cijm}
C^{(m)}_{ij} = \frac{[(N-\max\{ i,j \})m]_q [\min\{ i,j \}m]_q}{[Nm]_q}.
\end{equation}

The relations are:
\begin{align*}
\{ \s_i(z),\s_j(w) \} &= 2h(q-q^{-1}) C_{ij}\left( \frac{wq^{j-i}}{z}
\right) \s_i(z) \s_j(w) \\ &+ 2h \sum_{p=1}^i \delta\left(
\frac{w}{zq^{2p}} \right) \s_{i-p}(w) \s_{j+p}(z)
\\ &- 2h \sum_{p=1}^i \delta\left( \frac{wq^{2(j-i+p)}}{z} \right)
\s_{i-p}(z) \s_{j+p}(w),
\end{align*}
if $i\leq j$ and $i+j\leq N$; and
\begin{align*}
\{ \s_i(z),\s_j(w) \} &= 2h(q-q^{-1}) C_{ij}\left(
\frac{wq^{j-i}}{z} \right) \s_i(z) \s_j(w) \\  &+ 2h \sum_{p=1}^{N-j}
\delta\left( \frac{w}{zq^{2p}} \right) \s_{i-p}(w) \s_{j+p}(z)
\\ &- 2h \sum_{p=1}^{N-j} \delta\left( \frac{wq^{2(j-i+p)}}{z} \right)
\s_{i-p}(z) \s_{j+p}(w),
\end{align*}
if $i\leq j$ and $i+j>N$.

For example, in the case of $\sw_3$ we have the following relations:
\begin{align*}
\{ \s_1(z),\s_1(z) \} &= 2h (q-q^{-1}) \sum_{m\in\Z} \left( \frac{w}{z}
\right)^m \frac{[2m]_q [m]_q}{[3m]_q} \s_1(z) \s_1(w) \\ &+ 2h \delta \left(
\frac{w}{zq^2} \right) \s_2(z) - 2h \delta \left( \frac{wq^2}{z} \right)
\s_2(w),\\
\{ \s_1(z),\s_2(w) \} &= 2h (q-q^{-1}) \sum_{m\in\Z} \left( \frac{wq}{z}
\right)^m \frac{[m]_q^2}{[3m]_q} \s_1(z) \s_2(w) \\ &+ 2h \delta \left(
\frac{w}{zq^2} \right) - 2h \delta \left( \frac{wq^4}{z} \right),\\
\{ \s_2(z),\s_2(w) \} &= 2h (q-q^{-1}) \sum_{m\in\Z} \left( \frac{w}{z}
\right)^m \frac{[2m]_q [m]_q}{[3m]_q} \s_2(z) \s_2(w) \\ &+ 2h \delta \left(
\frac{w}{zq^2} \right) \s_1(w) - 2h \delta \left( \frac{wq^2}{z} \right)
\s_1(z).
\end{align*}

The Poisson algebra $\W_h(\sw_N)$ is a $q$--deformation of the classical
$\W$--algebra $\W(\sw_N)$. The asymptotic expansion of $\s_i(z)$ around
$h=0$ has the form $$\s_i(z) = \begin{pmatrix} N \\ i \end{pmatrix} + h^2
C_i W_2(z) + \ldots,$$ where $C_i$ is some coefficient, and $W_2(z)$ is the
quadratic (Virasoro) generating series of the classical ${\cal W}$--algebra
$\W(\sw_N)$. For each $n=2,\ldots,N-1$, one can find a combination of
$\s_i(z)$'s having expansion of the form $M + h^n W_n(z) + \ldots$, where
$M$ is a constant and $W_n(z)$ is a multiple of an $n$th order generating
series of $\W(\sw_N)$. The Poisson structure on $\W(\sw_N)$ can be
recovered from that on $\W_h(\sw_N)$.

\begin{rem}
If we replace $2h$ by $(q-q^{-1})$ in the formulas above, we will obtain a
Poisson algebra over the ring of rational functions in $q$.
\qed\end{rem}

\begin{rem} As in the case of $\sw_2$ (cf. \remref{dual}), the Poisson
algebra $\W_h(\sw_N)$ has another limit as $h \arr 0$, which can be
interpreted in terms of the central subalgebra of the algebra of functions
on the Poisson-Lie group dual to $\widehat{SL}_N$ \cite{next}.
\qed\end{rem}

\subsection{The center of $\wt{U}_h(\sun)$ at the critical level.}
Following \cite{RS}, cf. also \cite{FR,DE}, for any finite-dimensional
representation $W$ of $U_q(\G)$, one can construct matrices $L^\pm_W(z) =
(L^\pm_W(z))_{i,j=1,\ldots,\dim W}$ consisting of generating functions of
elements of $\wt{U}_h(\sun)\cri$. It is shown in \cite{RS} that the Fourier
coefficients of the power series
\begin{equation}    \label{gen}
\el^W(z) = \on{tr}_W q^{2\rho} L^+_W(z) L^-_W(zq^{-\k})^{-1}
\end{equation}
are central elements of $\wt{U}_h(\sun)\cri$.

In particular, for $\G=\sun$, let $\el_i(z) \equiv \el^{W_{\omega_i}}(z)$
be the generating function of central elements of $\wt{U}_h(\sun)\cri$
corresponding to the $i$th fundamental representation $W_{\omega_i}$ of
$U_q(\sun)$. Note that for $\G=\su$ we have $\el_1(z) = \el(zq^{-1})$,
where $\el(z)$ is given by \eqref{central}.

The next to leading term in the $h$--expansion of $\el_i(z)$
is a multiple of the Sugawara series of $\wt{U}(\sun)\cri$,
cf. \cite{DE}. Higher Sugawara elements of $\wt{U}(\sun)\cri$ can be
obtained from higher order terms of the expansions of $\el_i(z)$'s.

Let $Z_h(\sun)$ be the completion of the central subalgebra of
$\wt{U}_h(\sun)\cri$ generated by the coefficients of the series
$\el_1(z),\ldots,\el_{N-1}(z)$. Using the Wakimoto realization of
$\wt{U}_h(\sun)\cri$ \cite{AOS}, we obtain a homomorphism of Poisson
algebras $Z_h(\sun) \arr \HH$, which we call the $q$--deformation of the
Miura transformation.

Using the method of \remref{simple} we can find the image of $\el_1(z)$ in
$\HH$ by computing the matrix element of $\el_1(z)$ between the generating
vector of a Wakimoto module over $U_h(\sun)\cri$ and its dual. But for that
we only need the diagonal part of the ``Gauss decomposition'' of $L^\pm(z)$
\cite{DF} and a formula expressing the corresponding diagonal elements
$k_i^\pm(z), i=1,\ldots,N$, in terms of $\La_i(z), i=1,\ldots,N-1$. This
formula can be obtained from \cite{AOS} and \eqref{lai}. Explicit
computation shows that the image of $\el_1(z)$ in $\HH$ is equal to
$\s_1(z)$.

The generating functions $\el_i(z), i=2,\ldots,N-1$, corresponding to other
fundamental representations can be expressed in terms of $\el_1(z)$ by the
fusion procedure, cf. \cite{KRS}, \cite{CP} and references therein. Using
this procedure, we can show that the image of $\el_i(z)$ in $\HH$ is equal
to $\s_i(z)$ given by formula \eqref{qmap} for all $i=1,\ldots,N-1$. Thus,
we obtain that $Z_h(\sun)$ is isomorphic to $\W_h(\sw_N)$ as a Poisson
algebra. We will discuss this isomorphism in more detail in \cite{next}.

In conclusion of this section, recall that elements of the spectrum of
$\W(\sw_N)$ can be considered as $N$th order differential operators,
cf. \cite{DS}. The classical Miura transformation corresponds to splitting
of such an operator into a product of first order operators.

The spectrum of $\W_h(\sun)$ and the $q$--deformation of the Miura
transformation can be interpreted in a similar fashion. Namely, elements of
the spectrum of $\HH$ can be considered as first order $q$--difference
operators, and elements of the spectrum of $Z_h(\sun)$ can be considered as
$N$th order $q$--difference operators.

Indeed, let $Q_i(z), i=1,\ldots,N-1$, be solutions of the $q$--difference
equations $$({\cal D}_q - \La_i(z)) \frac{Q_i(z)}{Q_{i-1}(zq^{-2})} = 0,
\quad \quad i=1,\ldots,N,$$ where $[{\cal D}_q f](z) = f(zq^2)$ and we put
$Q_0(z)=Q_N(z)=1$. Then $Q_{N-1}(z)$ satisfies the $q$--difference equation
$$({\cal D}_q^N - \s_{N-1}(z) {\cal D}_q^{N-1} + \s_{N-2}(z) {\cal D}_q^{N-2}
- \ldots - (-1)^N \s_1 {\cal D}_q + (-1)^N) Q_{N-1}(zq^{-2}) = 0.$$

Thus, elements of the spectrum of $\W_h(\sw_N)$ can be considered as
$q$--difference operators of the form $${\cal D}_q^N - \s_{N-1}(z) {\cal
D}_q^{N-1} + \s_{N-2}(z) {\cal D}_q^{N-2} - \ldots - (-1)^N \s_1 {\cal D}_q +
(-1)^N.$$

\section{Generalization to other quantum affine algebras.}
Let us first adopt notation we used in the case of $\sun$ to the general
case. We pass to another set of generators of ${\cal H}_h(\sw_N)$, $y_i[n],
i=1,\ldots,N-1; n\in\Z$, such that $$\la_i[n] = q^{(i-1)n} y_i[n] - q^{in}
y_{i-1}[n], \quad \quad i=1,\ldots,N,$$ and $y_0[n]=y_N[n]=0.$ Using
formulas \eqref{pbl1} and \eqref{pbl2} we find the following Poisson
brackets:
\begin{equation}    \label{ypb}
\{ y_i[n],y_j[m] \} = 2h(q-q^{-1}) C_{ij}^{(m)} \delta_{n,-m},
\end{equation}
where $C^{(m)}_{ij}$ is given by \eqref{cijm},
\begin{equation}    \label{ya}
\{ y_i[n],a_j[m] \} = 2h [n]_q \delta_{n,-m} \delta_{ij}.
\end{equation}
Thus, the generators $y_i[n]$ are ``dual'' to the generators $a_i[n]$. In
fact, it is easy to see that
\begin{equation}    \label{inverse}
B^{(m)} C^{(m)} = [m]_q^2 I_{N-1},
\end{equation}
where $C^{(m)}=(C^{(m)}_{ij})_{i,j=1,\ldots,N-1}$ and
$B^{(m)}=([B_{ij}m]_q)_{i,j=1,\ldots,N-1}$, $(B_{ij})_{i,j=1,\ldots,N-1}$
being the Cartan matrix of $\sw_N$.

Introduce generating functions $$Y_i(z) = q^{-i(N-i)} \exp \left( -
\sum_{m\in\Z} y_i[m] z^{-m} \right).$$ We have: $$\La_i(z) = Y_i(zq^{-i+1})
Y_{i-1}(zq^{-i})^{-1}, \quad \quad i=1,\ldots,N,$$ where we put $Y_0(z) =
Y_N(z) = 1$. Note that $Y_i(z)$ can be written as
$Q_i(zq^{i+1})/Q_i(zq^{i-1})$ in terms of $Q_i(z)$ introduced at the end of
last section.

{}From formula \eqref{ypb} we find the Poisson brackets between $Y_i(z)$ and
$Y_j(w)$: $$\{ Y_i(z),Y_j(w) \} = 2h(q-q^{-1}) C_{ij} \left( \frac{w}{z}
\right) Y_i(z) Y_j(w),$$ where $C_{ij}(x)$ is given by formula \eqref{cij}.

We now define analogous generating functions $Y_i(z), i=1,\ldots,l$, for an
arbitrary simple Lie algebra $\g$. Namely, let $y_i[n], i=1,\ldots,l;
n\in\Z$, be the elements of ${\cal H}_h(\g)$ uniquely defined by
the Poisson bracket \eqref{ya}. We put $$Y_i(z) = q^{-2(\rho,\omega_i)}
\exp \left( - \sum_{m\in\Z} y_i[m] z^{-m} \right),$$ where
$\omega_i$ is the $i$th fundamental weight of $\g$. Then we have:
\begin{equation}    \label{yi}
\{ Y_i(z),Y_j(w) \} = 2h(q-q^{-1}) C_{ij} \left( \frac{w}{z} \right) Y_i(z)
Y_j(w),
\end{equation}
where $$C_{ij}(x) = \sum_{m\in\Z} C_{ij}^{(m)} x^m$$ and the matrix
$(C_{ij}^{(m)})_{i,j=1,\ldots,l}$ is defined by formula \eqref{inverse}
(with $N-1$ replaced by $l$) with respect to the symmetrized Cartan matrix
$B$ of $\g$.

\begin{rem} It is interesting that the functions $C_{ij}(x)$ appear in the
Thermodynamic Bethe Ansatz equations \cite{ORW,BR}.
\qed\end{rem}

For each dominant integral highest weight $\la$ of $\g$ there exists an
irreducible finite-dimensional representation $W_\la$ of $U_q(\G)$ which
satisfies the following property. Its restriction to the subalgebra
$U_q(\g)$ is completely reducible, the irreducible representation of
$U_q(\g)$ with highest weight $\la$ has multiplicity one, and all other
irreducible components of $W_\la$ have highest weights $\mu<\la$.

Let $\el_i(z), i=1,\ldots,l$, be the generating functions of central
elements of $\wt{U}_h(\G)\cri$ corresponding to $W_{\omega_i}$. Let
$Z_h(\G)$ be the central subalgebra of $\wt{U}_h(\G)\cri$ generated by the
coefficients of $\el_i(z), i=1,\ldots,l$. Recall from \S\S~1,2 that the
Miura transformation is the homomorphism of Poisson algebras from the
center of $\wt{U}(\G)$ to the Heisenberg-Poisson algebra ${\cal H}(\g)$.

\begin{conj}

\begin{enumerate}

\item[(a)]
{\em $Z_h(\G)$ is closed with respect to the Poisson structure on the
center of $\wt{U}_h(\G)\cri$.}

\item[(b)] {\em There exists a homomorphism of Poisson algebras $Z_h(\G)
\arr {\cal H}_h(\g)$, which is a deformation of the Miura transformation.}

\item[(c)]
{\em The formulas for the images $\s_i(z)$ of the generating functions
$\el_i(z)$ from $Z_h(\G)$ in ${\cal H}_h(\g)$ coincide with the formulas
for spectra of the corresponding transfer-matrices in integrable models
associated to $U_q(\G)$.}

\end{enumerate}
\end{conj}

Formulas for spectra of transfer-matrices in integrable models associated
to $U_q(\G)$ have been given in \cite{R1,R2,BR,KS} (although in different
normalizations). We will now describe $\s_i(z)$'s for all quantum affine
algebras of classical types via the series $Y_i(z), i=1,\ldots,l$ (we put
$Y_0(z)=1$). The Poisson brackets between $Y_i(z)$'s given by \eqref{yi}
uniquely determine the Poisson brackets between $\s_i(z)$'s.

\subsection{The series $\G=A_n^{(1)}$.}
Introduce
$$\La_i(z) = Y_i(zq^{-i+1}) Y_{i-1}(zq^{-i})^{-1}, \quad \quad
i=1,\ldots,n+1.$$

Let
$$\s_i(z) = \sum_{1\leq j_1 < \ldots < j_i\leq n+1} \La_{j_1}(z)
\La_{j_2}(zq^2) \ldots \La_{j_{i-1}}(zq^{2(i-2)}) \La_{j_i}(zq^{2(i-1)}),
\quad \quad i=1,\ldots,n.$$

\subsection{The series $\G=B_n^{(1)}$.}
Introduce
\begin{align*}
\La_i(z) &= Y_i(zq^{-i+1}) Y_{i-1}(zq^{-i})^{-1}, \quad \quad
i=1,\ldots,n-1, \\
\La_n(z) &= Y_n(zq^{-n+3/2}) Y_n(zq^{-n+1/2}) Y_{n-1}(zq^{-n})^{-1}, \\
\La_{n+1}(z) &= Y_n(zq^{-n+3/2}) Y_n(zq^{-n-1/2})^{-1}, \\
\La_{n+2}(z) &= Y_{n-1}(zq^{-n+1}) Y_n(zq^{-n+1/2})^{-1}
Y_n(zq^{-n-1/2})^{-1}, \\
\La_{2n-i+2}(z) &= Y_{i-1}(zq^{-2n+i+1}) Y_i(zq^{-2n+i})^{-1} \quad
\quad i=1,\ldots,n-1.
\end{align*}

Let
$$\s_i(z) = \sum_{\{ j_1,\ldots,j_i \} \in S} \La_{j_1}(z) \La_{j_2}(zq^2)
\ldots \La_{j_i}(zq^{2i-2}), \quad \quad i=1,\ldots,n-1,$$ where $S$ is the
set of $\{ j_1,\ldots,j_i \}$, such that $j_\al<j_{\al+1}$ or
$j_\al=j_{\al+1}=n+1, \al=1,\ldots,i-1$.

The formula for $\s_n(z)$, which corresponds to the spinor representation of
$B_n^{(1)}$, is more complicated:
$$\s_n(z) = \sum_{\si_1,\ldots,\si_n = \pm 1} b_{\si_1}(z|n)
b_{\si_2}(zq^{1-\si_1}|n-1) \ldots
b_{\si_n}(zq^{n-1-\si_1-\ldots-\si_{n-1}}|1),$$ where
\begin{align*}
b_1(z|1) &= Y_n(zq^{-n-1/2})^{-1}, \\ b_1(z|k) &= 1, \quad k=2,\ldots,n; \\
b_{-1}(z|1) &= Y_{n-1}(zq^{-n})^{-1} Y_n(zq^{-n+1/2}), \\
b_{-1}(z|k) &= Y_{n-k}(zq^{-n+k-1})^{-1} Y_{n+1-k}(zq^{-n+k}), \quad
k=2,\ldots,n.
\end{align*}

\subsection{The series $\G=C_n^{(1)}$.}
Introduce
\begin{align*}
\La_i(z) &= Y_i(zq^{-(i-1)/2}) Y_{i-1}(zq^{-i/2})^{-1}, \quad \quad
i=1,\ldots,n-1, \\
\La_n(z) &= Y_n(zq^{-(n-1)/2}) Y_{n-1}(zq^{-n/2})^{-1}, \\
\La_{n+1}(z) &= Y_{n-1}(zq^{-(n+2)/2}) Y_n(zq^{-(n+3)/2})^{-1}, \\
\La_{2n-i+1}(z) &= Y_{i-1}(zq^{-(2n-i+2)/2}) Y_i(zq^{-(2n-i+3)/2})^{-1},
\quad \quad i=1,\ldots,n-1.
\end{align*}

Let
$$\s_i(z) = \sum_{\{ j_1,\ldots,j_i \} \in S} \La_{j_1}(z) \La_{j_2}(zq)
\ldots \La_{j_i}(zq^{i-1}), \quad \quad i=1,\ldots,n,$$ where $S$ is the
set of $\{ j_1,\ldots,j_i \}$, such that $1\leq j_1 < \ldots < j_i\leq 2n$
and if $j_\al = l, j_\beta = 2n+1-l$ for some $l=1,\ldots,n$, then $l\leq
n+\al-\beta$.

\subsection{The series $\G=D_n^{(1)}$.}
Introduce
\begin{align*}
\La_i(z) &= Y_i(zq^{-i+1}) Y_{i-1}(zq^{-i})^{-1}, \quad \quad
i=1,\ldots,n-2, \\
\La_{n-1}(z) &= Y_n(zq^{-n+2}) Y_{n-1}(zq^{-n+2}) Y_{n-2}(zq^{-n+1})^{-1},
\\
\La_n(z) &= Y_{n-1}(zq^{-n+2}) Y_n(zq^{-n})^{-1}, \\
\La_{n+1}(z) &= Y_n(zq^{-n+2}) Y_{n-1}(zq^{-n})^{-1}, \\
\La_{n+2}(z) &= Y_{n-2}(zq^{-n+1}) Y_{n-1}(zq^{-n})^{-1}
Y_n(zq^{-n})^{-1}, \\
\La_{2n-i+1}(z) &= Y_{i-1}(zq^{-2n+i+2}) Y_i(zq^{-2n+i+1})^{-1}.
\end{align*}

Let
$$\s_i(z) = \sum_{\{ j_1,\ldots,j_i \} \in S} \La_{j_1}(z) \La_{j_2}(zq^2)
\ldots \La_{j_i}(zq^{2i-2}), \quad \quad i=1,\ldots,n-2,$$ where $S$ is the
set of $\{ j_1,\ldots,j_i \}$, such that $j_\al<j_{\al+1}$ or
$j_\al=j_{\al+1}+1=n+1, \al=1,\ldots,i-1$.

The formulas for $\s_{n-1}(z)$ and $\s_n(z)$, which correspond to the spinor
representations of $D_n^{(1)}$, are more complicated. In these formulas the
subscript $\ep$ means $n$, if $\ep=+$, and $n-1$, if $\ep=-$. Thus,
$\s_+(z)=\s_n(z), \s_-(z)=\s_{n-1}(z)$. Now let
$$\s_\ep(z) = \sum_{\si_1,\ldots,\si_{n-1}=\pm 1} b_{\si_1}^\ep(z|n)
b_{\si_2}^{\ep\si_1}(zq^{1-\si_1}|n-1) \ldots b_{\si_{n-1}}^{\ep
\si_1\ldots\si_{n-1}}(zq^{n-2-\si_1-\ldots-\si_{n-2}}|2),$$ where
\begin{align*}
b^\ep_1(z|2) &= Y_\ep(zq^{-n})^{-1}, \\ b^\ep_1(z|k) &= 1, \quad
k=3,\ldots,n; \\
b^\ep_{-1}(z|2) &= Y_{n-2}(zq^{-n+1})^{-1} Y_\ep(zq^{-n+2}), \\
b^\ep_{-1}(z|k) &= Y_{n-k}(zq^{-n+k-1})^{-1} Y_{n+1-k}(zq^{-n+k}), \quad
k=3,\ldots,n.
\end{align*}

\end{document}